\documentclass{amsart}
\usepackage{amssymb}
\usepackage{pictexwd}

\newtheorem{theorem}{Theorem}[section]
\newtheorem{proposition}[theorem]{Proposition}
\newtheorem{corollary}[theorem]{Corollary}
\newtheorem{lemma}[theorem]{Lemma}

\newtheorem{preremark}[theorem]{Remark}
\newtheorem{predefinition}[theorem]{Definition}
\newtheorem{preexample}[theorem]{Example}
\newenvironment{remark}{\begin{preremark}\rm}{\end{preremark}}
\newenvironment{definition}{\begin{predefinition}\rm}{\end{predefinition}}
\newenvironment{example}{\begin{preexample}\rm}{\end{preexample}}

\newcounter{itemcounter}
\newenvironment{items}{
    \begin{list}{\alph{itemcounter})}
    {\usecounter{itemcounter}\setlength{\topsep}{3 pt}
    \setlength{\leftmargin}{2 em}
    \setlength{\partopsep}{0 pt}\setlength{\itemsep}{0 pt}
       \setlength{\labelwidth}{2 em}
    }}{\end{list}}

\let\phi=\varphi
\let\rho=\varrho
\let\theta=\vartheta
\let\epsilon=\varepsilon

\def\LT{\mathop{\rm LT}\nolimits}

\def\Mat{\mathop{\rm Mat}\nolimits}

\def\Supp{\mathop{\rm Supp}\nolimits}
\def\Spec{\mathop{\rm Spec}\nolimits}
\def\Proj{\mathop{\rm Proj}\nolimits}
\def\Proj{\mathop{\rm Proj}\nolimits}

\def\calA{\mathcal{A}}
\def\calO{\mathcal{O}}
\def\calZ{\mathcal{Z}}
\def\calS{\mathcal{S}}

\def\bbA{\mathbb{A}}
\def\bbB{\mathbb{B}}
\def\bbG{\mathbb{G}}

\def\bbT{\mathbb{T}}
\def\bbN{\mathbb{N}}
\def\bbP{\mathbb{P}}

\def\mapdown#1{\Big\downarrow\rlap{$\vcenter{\hbox{$\scriptstyle
    #1$}}$}}
\def\mapup#1{\Big\uparrow\rlap{$\vcenter{\hbox{$\scriptstyle
    #1$}}$}}

\let\To=\longrightarrow
\def\TTo#1{\mathop{\longrightarrow}\limits ^{#1}}
\def\iso#1{\mathop{\simeq}\limits ^{^{#1}}}

\def\mapdown#1{\Big\downarrow\rlap{$\vcenter{\hbox{$\scriptstyle
    #1$}}$}}

\def\mapup#1{\Big\uparrow\rlap{$\vcenter{\hbox{$\scriptstyle
    #1$}}$}}

\let\goth=\mathfrak

\def\cocoa{\mbox{\rm
    C\kern-.13em o\kern-.07 em C\kern-.13em o\kern-.15em A}}

\begin{document}

\title{On Border Basis and   Gr\"obner Basis Schemes}


\author{L. Robbiano}
\address{Dipartimento di Matematica, Universit\`a di Genova,
\newline \phantom{\;\;\;\,} Via Dodecaneso 35,
I-16146 Genova, Italy}
\email{robbiano@dima.unige.it}

\date{\today}

\begin{abstract}
Hilbert schemes of zero-dimensional ideals in a polynomial ring 
can be covered with suitable affine open subschemes whose 
construction is achieved using border bases.
Moreover, border bases have proved to be an excellent tool for
describing zero-dimensional ideals when the coefficients are inexact.
And in this situation they show a clear advantage with respect to 
Gr\"obner bases which, nevertheless, can  also be used
in the study of Hilbert schemes, since they provide tools for
constructing suitable stratifications. 

In this paper we compare
Gr\"obner basis schemes with border basis schemes.
It is shown that
Gr\"obner basis schemes and their associated universal families 
can be viewed as weighted projective schemes.
A first consequence of our approach is 
the proof that  {\it all}\/ the ideals which 
define a Gr\"obner basis scheme and are
obtained using Buchberger's Algorithm, are equal. 
Another result is  that if the origin (i.e.\ the point 
corresponding to the unique monomial ideal) in the 
Gr\"obner basis scheme  is smooth, then the scheme 
itself is isomorphic to an affine space. 
This fact represents a remarkable difference between 
border basis and Gr\"obner basis schemes. 
Since it is  natural to look for situations 
where a Gr\"obner basis scheme and the corresponding border 
basis scheme are equal, we address the issue,
provide an answer, and exhibit some consequences.
Open problems are discussed at the end of the paper.

\end{abstract}


\maketitle


\section{Introduction}

This paper has three main sources and one ancestor. 
Let me be more specific.

{\it Source 1.}\/ Given a zero-dimensional  ideal $I$ in a polynomial ring, 
and assuming that the coefficients of the 
generating polynomials are inexact, what is the 
best way of describing $I$?
The idea that Gr\"obner bases are not suitable for 
computations with inexact data
has been brought to light by Stetter (see~\cite{S}) 
and other numerical analysts.
 Gr\"obner bases are inadequate
due to the rigidity imposed by the term ordering. 
The class of border bases is more promising.
A~pioneering paper on border bases is~\cite{M} and
a~detailed description  is contained in Section 6.4 of~\cite{KR2}.

{\it Source 2.}\/ The possibility of parametrizing families 
of schemes with a scheme is a
remarkable peculiarity of algebraic geometry.
Hilbert schemes are one instance of this phenomenon, 
and consequently are widely studied.  If we let $P = K[x_1, \dots, x_n]$, 
Hilbert schemes of zero-dimensional ideals in $P$ can be covered 
by affine open subschemes which parametrize
all the subschemes $\Spec(P/I)$ of the affine space~$\bbA^n_K$ 
such that $P/I$ has a fixed basis.  
What is  interesting is that the construction of such subschemes
is performed using border bases 
(see for instance~ \cite{Hu1}, \cite{Hu2}, and~ \cite{MS}).

{\it Source 3.}\/  Despite their inability to treat inexact 
data well, Gr\"obner bases can  nevertheless be used
in the study of Hilbert schemes, since with their help it is possible to
construct suitable stratifications. 
Among the vast literature on this subject let me mention
the two fairly recent articles~\cite{CV} and~\cite{NS} 
and the bibliography quoted therein.

The three main sources are now described; 
it remains to reveal the ancestor.
It is paper~\cite{KR3} where we tried to extend 
to border bases a very nice property of 
Gr\"obner bases, the possibility of connecting 
every ideal to its leading term ideal via a flat
deformation. We were able to get partial results, so that
the connectedness of  border basis schemes is 
still an open problem (see Question 2 at the end of this paper). 

So what is the content of the next pages? The main idea is to compare 
Gr\"obner basis schemes (see Definition~\ref{genericGB})  
with border basis schemes. We also define a universal family 
(see Definition~\ref{universalGB}), and the first main result is 
Theorem~\ref {homogeneousfamily}   where it is shown that 
Gr\"obner basis schemes and their associated universal families 
can be endowed with a graded structure where the indeterminates 
have positive weights. In other words, they can be viewed as 
weighted projective schemes (see Remark~\ref{specandproj}). 
The second main result is Theorem~\ref{GandB} where the 
comparison of the two schemes is fully described. In particular,
it is shown that Gr\"obner basis schemes can be obtained 
as sections of border basis schemes with suitable linear spaces.
Since our description of  Gr\"obner basis schemes is not directly 
linked to the concept of Gr\"obner basis, we prove in
Corollary~\ref{trueGfamily} that indeed our definition 
is well-placed.

Section 3 is devoted to exhibiting some consequences of 
the above mentioned results. Let me explain the first one. 
In the literature Gr\"obner basis schemes are mostly described
using Buchberger's Algorithm. 
However, this approach has a drawback, since the 
reduction process in the algorithm is far from being unique, 
and the consequence is that the description
of the Gr\"obner basis scheme is {\it a priori}\/ not canonical.  
A first consequence of our approach is 
the proof that  {\it all}\/ the ideals obtained using 
Buchberger's Algorithm are 
equal (see Proposition~\ref{democracy}) 
and coincide with the  ideal defined in this paper
(see Proposition~\ref{BBisreduction}).

Another remark is made in Corollary~\ref{affinecell} 
where it is shown that if the origin (i.e.\ the point 
corresponding to the unique monomial ideal) in the 
Gr\"obner basis scheme  is smooth, then the scheme 
itself is isomorphic to an affine space. 
This fact represents a remarkable difference between 
border basis and Gr\"obner basis schemes (see Example~\ref{square}). 

After Theorem~\ref{GandB} it is  natural to look for situations 
where a Gr\"obner basis scheme and the corresponding border 
basis scheme are the same. The answer is given in 
Proposition~\ref{iso} and a nice consequence is shown 
in Corollary~\ref{stick}.

\smallskip
Doing mathematics is looking for solutions to problems, 
a process which inevitably sparks new questions. 
This paper is no exception; in particular, two open 
questions are presented at the end of Section 3.

\begin{flushright}
\small\it
Judge others by their questions\\
 rather than by their answers.\\
\rm (Fran\c cois-Marie Arouet (Voltaire))
\end{flushright} 

\medskip

Unless explicitly stated otherwise, we use definitions and notation
introduced in~\cite{KR1},~\cite{KR2},~\cite{KR3}. All the experimental computation
was done with the computer algebra system \cocoa\  (see~\cite{CoCoA}).

\bigskip

\section{Border Basis and Gr\"obner Basis Schemes}
\label{BBand GB}

In the following we let $K$ be a field, $P=K[x_1,\dots,x_n]$ a
polynomial ring, and~$I\subset P$ a zero-dimensional ideal.
Recall that an {\em  order ideal}~$\calO$
is a finite set of terms in
$\bbT^n= \bbT(x_1,\dots,x_n)=
\{x_1^{\alpha_1}\cdots x_n^{\alpha_n} \mid \alpha_i\ge 0\}$
such that all divisors of a term in~$\calO$
are also contained in~$\calO$.
The set 
$\partial\calO=(x_1\calO\cup \cdots \cup x_n\calO) \setminus \calO$
is called the {\em border} of~$\calO$.

\begin{definition}
Let $\calO=\{t_1,\dots,t_\mu\}$ be an order ideal and 
$\partial\calO=\{b_1,\dots,b_\nu\}$ its border.

\begin{items}
\item A set of polynomials $\{g_1,\dots,g_\nu\}\subseteq I$
is called an {\em $\calO$-border prebasis} of~$I$ if it is of
the form
$g_j=b_j-\sum_{i=1}^\mu a_{ij}t_i$ with $a_{ij}\in K$.

\goodbreak
\item An $\calO$-border prebasis of~$I$ is called an
{\em $\calO$-border basis} of~$I$ if
$P=I\oplus \langle \calO \rangle_K$.

\end{items}
\end{definition}

It is known that if $I$ has an  $\calO$-border basis,
then such $\calO$-border basis of $I$ is
unique (see~\cite{KR2} Proposition 6.4.17).

\begin{proposition}{\bf (Border Bases and Multiplication Matrices)}
\label{formalmult}\\
Let $\calO=\{t_1,\dots,t_\mu\}$  be an order ideal of monomials,
let the set $\{g_1,\dots,g_\nu\}$ be
an $\calO$-border prebasis, and let $I $ be the ideal
generated by $\{g_1,\dots,g_\nu\}$. Then the following conditions are
equivalent
\begin{items}
\item The set $\{g_1,\dots,g_\nu\}$ is the $\calO$-border basis
of $I$.

\item The formal multiplication matrices of $\{g_1,\dots,g_\nu\}$
are pairwise commuting.
\end{items}
\end{proposition}

\proof See~\cite{KR2}, Definition 6.4.29 and  Theorem 6.4.30.
\endproof
\goodbreak

\begin{definition}\label{defBBS}
Let $\{c_{ij} \mid 1\le i\le \mu,\;
1\le j\le\nu\}$ be a set of new indeterminates.

\begin{items}
\item The {\em generic $\calO$-border prebasis}
is the set of polynomials
$G=\{g_1,\dots,g_\nu\}$ in~$K[x_1,\dots,x_n,c_{11},\dots,c_{\mu\nu}]$
given by
$$
g_j = b_j -\sum_{i=1}^\mu c_{ij}t_i
$$

\item For $k=1,\dots,n$, let $\calA_k \in\Mat_{\mu}(K[c_{ij}])$
be the
$k^{\rm th}$ formal multiplication matrix associated to~$G$
(cf.~\cite{KR2}, Definition 6.4.29).
It is also called the $k^{\rm th}$ {\em generic
multiplication matrix}\/ with respect to~$\calO$.

\item The ideal of 
$K[c_{11}, \dots, c_{\mu \nu}]$ generated by the entries of 
$\calA_k \calA_\ell -\calA_\ell \calA_k$ with
$1\le k<\ell\le n$ defines an affine subscheme of
$\bbA^{\mu\nu}$
which will be denoted by~$\bbB_{\calO}$ and called 
the {\em $\calO$-border basis scheme}. 
Its defining ideal will
be denoted by $I(\bbB_{\calO})$,
and its coordinate ring
$K[c_{11},\dots,c_{\mu\nu}]/I(\bbB_\calO)$
will be denoted by~$B_{\calO}$.

\end{items}
\end{definition}

The reason why it is called the $\calO$-border 
basis scheme is the following.
When we apply the substitution $\Sigma(c_{ij})=\alpha_{ij}$
to~$G$, a point $(\alpha_{ij})\in K^{\mu\nu}$ yields a
border basis   if and only if 
$\Sigma(\calA_k)\,\Sigma(\calA_\ell)=
\Sigma(\calA_\ell)\,\sigma(\calA_k)$
for $1\le k<\ell\le n$ (see~Proposition~\ref{formalmult}). Thus the
$K$-rational points of~$\mathbb
B_{\calO}$ are in 1--1 correspondence with the
$\calO$-border bases of
zero-dimensional ideals in~$P$,
and therefore are in 1--1 correspondence with 
{\it all zero-dimensional ideals $I$ in $P$
such that $\overline{\calO}$ is a basis of 
$P/I$ as a~$K\!$-vector space}.

\bigskip
Next, we are going to define
$(\calO, \sigma)$-Gr\"obner basis schemes, and
to do this an extra bit of notation is required.
Let $\calO=\{t_1,\dots,t_\mu\}$ be an order ideal.
Then the set of  minimal generators of the
monoideal~$\mathbb T^n\setminus \calO$ (also called
the {\em corners}\/ of~$\calO$) is denoted
by~$c\calO$, and
we denote by~$\eta$ the cardinality of~$c\calO$.
Since~$c\calO \subseteq \partial\calO$,
it follows that $\eta \le \nu$, and
we label the elements in $\partial\calO$ so that
$c\calO=\{b_1, \dots, b_\eta\}$.

We let $\sigma$ be a term ordering on $\mathbb T^n$
and recall that if $I$ is an ideal in the polynomial ring $P$, we denote
the order ideal
$\mathbb{T}^n\setminus \LT_\sigma(I)$ by
$\calO_\sigma(I)$. Moreover, 
we denote by $S_{\calO, \sigma}$ the set
$\{c_{ij} \in \{c_{11}, \dots, c_{\mu\nu} \}   \ | \  b_j >_\sigma t_i   
\}$,
by~$L_{\calO, \sigma}$ the ideal generated
by~$\{c_{11},\dots, c_{\mu\nu}\}\setminus S_{\calO, \sigma}$
in~$K[c_{11},\dots,c_{\mu\nu}]$,
by $S_{c\calO, \sigma}$ the intersection
$S_{\calO, \sigma}\cap \{c_{11},\dots,c_{\mu\eta}\}$,
and by $L_{c\calO, \sigma}$ the
ideal generated by~$\{c_{11}, \dots, c_{\mu\eta}\}\setminus S_{c\calO, \sigma}$
in~$K[c_{11},\dots,c_{\mu\eta}]$.
Furthermore we denote  the
cardinality of $S_{c\calO,\sigma}$ by  $s(c\calO,\sigma)$.


\goodbreak

\begin{definition}\label{genericGB}
For $j = 1,\dots, \nu$ we define $g_j^*$ in the following way.
$$
g_j^* \ =\  b_j -\sum_{ \{ i \ | \ b_j >_{_\sigma} t_i \} }c_{ij}t_i \ = \ 
b_j -\sum_{ c_{ij} \in   S_{\calO,\sigma}\cap \{c_{1j}, \dots, c_{\mu j} \} }  c_{ij}t_i
$$

\begin{items}
\item 
The {\em generic $(\calO, \sigma)$-Gr\"obner prebasis} is
the set of polynomials $\{g_1^*, \dots, g_\eta^*\}$.

\item
The ideal 
$\big(L_{\calO, \sigma}+ I(\bbB_{\calO})\big)\cap K[S_{c\calO,\sigma}]$
of  $K[S_{c\calO,\sigma}]$
defines an affine subscheme of
$\bbA^{s(c\calO,\sigma)}$
which will be denoted by $\bbG_{\calO, \sigma}$
and called  the{\em $(\calO, \sigma)$-Gr\"obner basis scheme}.
The defining ideal 
$ \big(L_{\calO, \sigma}+ I(\bbB_{\calO})\big)\cap K[S_{c\calO,\sigma}]$ will
be denoted by~$I(\bbG_{\calO, \sigma})$ and the coordinate ring
$K[S_{c\calO,\sigma}]/I(\bbG_{\calO, \sigma})$
will be denoted by~$G_{\calO, \sigma}$.

\end{items}
We observe that $g_j^*$ is obtained from $g_j$ by 
setting to zero all the indeterminates in 
$L_{\calO, \sigma}\cap \{c_{1j}, \dots, c_{\mu j}\}  $.

\end{definition}

\goodbreak

\begin{example}
We examine the inclusion $c\calO \subseteq \partial\calO$.
If $\calO = \{1,x,y,xy\}$ then $c\calO = \{x^2, y^2\}$ while
$\partial\calO = \{x^2, y^2, x^2y, xy^2 \}$,
so that $c\calO \subset \partial\calO$.
On the other hand, if $\calO = \{1,x,y\}$ then
$c\calO = \partial\calO =\{x^2, xy, y^2\}$.

Returning to $\calO = \{1,x,y,xy\}$
we observe that $t_1 = 1$, $t_2 = x$, $t_3=y$, $t_4 = xy$,
$b_1 = x^2$, $b_2 = y^2$, $b_3 = x^2y$, $b_4 = xy^2$.
Let $\sigma = {\tt DegRevLex}$, so that  $x>_\sigma y$. 
Then $L_{\calO, \sigma} = L_{c\calO, \sigma} = (c_{42} )$, $g_1^* = g_1$, 
$g_2^* =y^2-( c_{12} + c_{22}x+c_{32}y)$, $g_3^*=g_3$, $g_4^*=g_4$.

\end{example}

\bigskip

Having introduced the Gr\"obner basis scheme, we  define
a naturally associated universal family. To  this end
we recall the following definition taken from~\cite{KR3} and
extend it.

\begin{definition}\label{universalGB}
The ring
$K[x_1,\dots,x_n, c_{11},\dots,c_{\mu\nu}]/
\big(I(\bbB_{\calO})+(g_1,\dots, g_\nu )\big)$
will be denoted by~$U_{\calO}$. The ring
$K[x_1,\dots,x_n, S_{c\calO,\sigma}]/
\big(I(\bbG_{\calO, \sigma})+(g_1^*,\dots, g_\eta^* )\big)$
will be denoted by~$U_{\calO, \sigma}$.

\begin{items}
\item
The natural homomorphism of $K$-algebras
$\Phi:\; B_{\calO} \To  U_{\calO}$
is called the {\it universal  $\calO
$-border basis family}.

\item The natural homomorphism of $K$-algebras
$\Psi:\; G_{\calO,\sigma} \To  U_{\calO,\sigma}$
is called the {\it universal  $(\calO,\sigma)$-Gr\"obner basis family}.

\item
The induced homomorphism of $K$-algebras
$B_{\calO}/ \overline{L}_{\calO,\sigma} \;\longrightarrow\;
U_{\calO} /\overline{L}_{\calO,\sigma}$ 
will be denoted by $\overline{\Phi}$.
\end{items}

\end{definition}

\begin{remark}\label{posweights}
It is known (see~\cite{Bay}, and ~\cite{E} Exercise 15.12, p. 370) that 
given power products  $t, t_1, \dots, t_r \in \bbT^n$ and 
a term ordering $\sigma$ such that 
$t>_\sigma t_i$ for $i =1, \dots r$, then there exists a 
system~$V$ of positive weights on~$x_1, \dots, x_n$ (i.e\ a matrix
$V \in \Mat_{1,n}(\bbN_+)$) such that $\deg_V(t) > \deg_V(t_i)$ 
for $i = 1, \dots, r$.
\end{remark}

We are ready to prove an important property of some ideals 
described before. To help the reader, we observe that for simplicity we write 
$\mathbf{x}$ for  $x_1, \dots, x_n$ and~$\mathbf{c}$ 
for $c_{11}, \dots c_{\mu\nu}$. 

\begin{theorem}\label{homogeneousfamily}  
\label{iboHom}
There exist a system~$W$  of positive weights on the elements of
 $S_{c\calO, \sigma}$,
a system~$\overline{W}$ of positive weights on the elements of
$S_{\calO, \sigma}$, and a system~$V$ of positive weights on
$\mathbf{x}$
such that the following conditions hold true.

\begin{items}
\item The system $\overline{W}$ is an extension of the system~$W\!$.

\item The ideal $I(\bbG_{\calO, \sigma})$
in $K[S_{c\calO, \sigma}]$ is $W$-homogeneous.

\item The ideal $I(\bbG_{\calO, \sigma}) + (g_1^*, \dots, g_\eta^*)$
in $K[\mathbf{x}, S_{c\calO, \sigma}]$ 
is $(V,W)$-homogeneous.

\item The image of $I(\bbB_{\calO})$
in   $K[S_{\calO, \sigma}] \cong 
K[\mathbf{c}]/L_{\calO, \sigma}$
is  $\overline{W}$-homogeneous.

\item The image of $I(\bbB_{\calO})+(g_1^*, \dots, g_\nu^*)$
in  $K[\mathbf{x}, S_{\calO, \sigma}] \cong 
K[\mathbf{x}, \mathbf{c}]/L_{\calO, \sigma}$
is  $(V,\overline{W})$-homogeneous.

\end{items}
\end{theorem}

\proof
The definition of $S_{c\calO, \sigma}$ and 
Remark~\ref{posweights} imply that there exists a system $V$ 
of positive weights on  $\mathbf{x}$
such  that $\deg_V(b_j) > \deg_V(t_{i})$
for every $j =1,\dots,\eta$ and 
every $t_i \in \Supp(g_j^*-b_j)$.
We define~$W$ by giving the
$c_{ij}\!$'s suitable positive weights, so that
all  elements $g_j^*$ in the generic
$(\calO, \sigma)$-Gr\"obner prebasis are
$(V, W)$-homogeneous
when they are viewed as polynomials in
$K[\mathbf{x}, S_{c\calO, \sigma}]$.

Then we choose a $\deg_{(V,W)}$-compatible
term ordering $\overline{\sigma}$ on
$\mathbb T(\mathbf{x}, S_{c\calO, \sigma})$
with the property that for every 
$t, t' \in \mathbb T(S_{c\calO, \sigma})$,
$x_1^{a_1}\cdots x_n^{a_n}\,t >_{\overline{\sigma}}
x_1^{b_1}\cdots x_n^{b_n}\,t'$ if
they have the same $(V,W)$-degree and
$x_1^{a_1}\cdots x_n^{a_n}  >_\sigma x_1^{b_1}\cdots x_n^{b_n}$.
If we use the $\overline{\sigma}$-division algorithm with respect to the
tuple
$(g_1^*,  \dots, g_\eta^*)$,
we can express every element
$b_j\in \partial \calO \setminus c \calO$
as a linear combination of those  elements in~$\calO$
which are $\sigma$-smaller than $b_j$.
Since all the~$g_i^*$ are monic and homogeneous,
the coefficients~$h_{ij}$ of these linear combinations are
homogeneous polynomials  in the~$c_{ij}\!$'s.
We define $\overline{W}$ by putting
$\deg_{\overline{W}}(c_{ij}) = \deg_W(h_{ij})$ for 
$c_{ij} \in S_{\calO, \sigma} \cap \{c_{1j}, \dots, c_{\mu j} \} $
and $j=\eta+1,\dots,\nu$.
We observe that $\overline{W}$ does not depend
on the choice of the order in the division algorithm,
it only depends on  $\calO$, $\sigma$, $V\!$.
At this point we have proved statement $a)$ and have shown that 
the polynomials $g_1^*, \dots, g_\nu^*$ are 
$(V, \overline{W})$-homogeneous which implies
that $d)$ and $e)$ are equivalent.
Moreover, we observe that~$b)$ follows from $d)$, while 
$c)$ and $d)$ follow from $e)$, 
so we only need to prove~$d)$.  Multiplication by $x_i$ yields a graded
homomorphism between $(V, \overline{W})$-graded free
$K[\mathbf{x}, \mathbf{c}]/L_{\calO, \sigma}$-modules,
therefore the multiplication matrices are homogeneous
(see~\cite{KR2}, Definition 4.7.1 and Proposition 4.7.4).
Consequently, the image of the ideal  $I(\bbB_{\calO})$
modulo $L_{\calO,\sigma}$ is $\overline{W}$-homogeneous and the
proof is complete.
\endproof

In the sequel we consider the following commutative diagram
of canonical homomorphisms

$$
\begin{matrix}
G_{\calO,\sigma} &\TTo{\phi} & B_{\calO}/\overline{L}_{\calO,\sigma}     \cr
&&& \cr 
\mapdown{\Psi} && \mapdown{\overline{\Phi}}  \cr
&&\cr
U_{\calO,\sigma}\ &\TTo{\theta}  &
U_{\calO}/\overline{L}_{\calO,\sigma} 
\end{matrix}
\eqno (1)
$$

i.e.\ 

$$
\begin{matrix}
K[S_{c\calO,\sigma}]/I(\bbG_{\calO, \sigma}) 
&\TTo{\phi} &
K[\mathbf{c}]/\big( L_{\calO, \sigma} + I( \bbB_{\calO})  \big)     \cr
&&\cr
\mapdown{\Psi} && \mapdown{\overline{\Phi}} \cr
&&\cr
K[\mathbf{x}, S_{c\calO, \sigma}]/\big( I(\bbG_{\calO, \sigma})     
+(g_1^*, \dots, g_\eta^*)\big)
&
\TTo{\theta}  &
K[\mathbf{x}, \mathbf{c}] / 
\big(L_{\calO,\sigma} + I(\bbB_{\calO}) + (g_1, \dots g_\nu)  \big)\cr
\end{matrix}
$$

\medskip
\noindent

We recall the equality
$I(\bbG_\calO) =
\big( L_{\calO, \sigma} + I(\bbB_{\calO}) \big) \cap K[S_{c\calO,\sigma}] $
from which the homomorphism $\phi$ derives.
The homomorphism $\theta$ is obtained as follows: 
let 
${\Theta: K[\mathbf{x}, S_{c\calO, \sigma}] \To K[\mathbf{x}, \mathbf{c}]}$ 
be the natural inclusion of polynomial rings. Then clearly
$I(\bbG_{\calO, \sigma})     +(g_1^*, \dots, g_\eta^*)
\subseteq 
\Theta^{-1}\big(L_{\calO,\sigma} + I(\bbB_{\calO}) + (g_1, \dots g_\nu)  \big) $.

\medskip 
We are ready to state the main result of this section. 
To prove it we are going to make extensive use of the above 
diagram $(1)$.

\begin{theorem}{\bf (Gr\"obner and Border)}
\label{GandB}\\
Let $\calO=\{t_1,\dots,t_\mu\}$  be an order ideal of monomials in $P$ and
let $\sigma$ be a term ordering on $\bbT^n$.

\begin{items}
\item The classes of the elements in $\calO$ form a 
$B_\calO/\overline{L}_{\calO,\sigma}$-module basis 
of~$U_{\calO}/\overline{L}_{\calO,\sigma}$. 

\item The classes of the elements in $\calO$ form a 
$G_{\calO, \sigma}$-module basis 
of~$U_{\calO, \sigma}$. 

\item We have the equality
$I(\bbG_{\calO, \sigma})+(g_1^*,\dots, g_\eta^*) =
\theta^{-1}\big(L_{\calO,\sigma} + I(\bbB_{\calO}) + 
(g_1, \dots g_\nu)  \big)  $.

\item The maps $\phi$ and $\theta$ in the above diagram 
are isomorphisms.
\end{items}

\end{theorem}

\proof
We observe that $\phi$ is injective by definition.  The fact that 
${\Phi: B_\calO \To U_\calO}$ is free with basis $\overline{\calO}$
and injective is proved in~\cite{KR3}, Theorem 3.4.
Passing to the quotient modulo $\overline{L}_{\calO, \sigma}$, we deduce that 
$\overline{\Phi} : B_\calO/\overline{L}_{\calO, \sigma}: 
\To U_\calO/\overline{L}_{\calO, \sigma}$ is free with 
basis $\overline{\calO}$ and injective, so that a) is proved.
We divide the proof of b)  into two claims.

{\it Claim 1. $\overline{\calO}$ generates}. 
According to Theorem~\ref{homogeneousfamily},
we may choose positive weights $W$ on the elements of $S_{c\calO, \sigma}$ 
and positive weights $V$ on  $\mathbf{x}$ so that 
$I(\bbG_{\calO,\sigma}) $ is a~$W$-homogeneous ideal 
of $K[S_{c\calO,\sigma}]$ and 
$I(\bbG_{\calO,\sigma}) + (g_1^*, \dots, g_\eta^*)$ is a 
$(V,W)$-homogeneous ideal of $K[\mathbf{x}, S_{c\calO,\sigma}]$. 
Following the lines of the proof of 
Theorem~\ref{homogeneousfamily}, we choose
a $\deg_{(U,W)}$-compatible
term ordering $\overline{\sigma}$ on
$\mathbb T(\mathbf{x}, S_{c\calO, \sigma})$
with the property that for every
$t, t' \in \mathbb T(S_{c\calO, \sigma})$,
$x_1^{a_1}\cdots x_n^{a_n}\,t >_{\overline{\sigma}}
x_1^{b_1}\cdots x_n^{b_n}\,t'$ if
they have the same $(V,W)$-degree and
$x_1^{a_1}\cdots x_n^{a_n}  >_\sigma x_1^{b_1}\cdots x_n^{b_n}$.
We observe that $\calO$ is the complement in $\bbT^n$ 
of the monoideal generated by $c\calO$, hence
if we use the $\overline{\sigma}$-division algorithm with respect to the
tuple $(g_1^*,  \dots, g_\eta^*)$,
we can express every polynomial in $K[\mathbf x, S_{c\calO, \sigma}]$
as a linear combination of  elements in~$\calO$, 
modulo $(g_1^*,  \dots, g_\eta^*)$.

{\it Claim 2. $\overline{\calO}^{\mathstrut}$ is linearly 
independent over $G_{\calO, \sigma}$}.
Let $f = \sum_{i=1}^\mu f_it_i \in K[\mathbf x, S_{c\calO, \sigma}]$ 
and assume that $f = 0 \  {\rm modulo} \ 
\big(I(\bbG_{\calO,\sigma}) +(g^*_1, \dots, g_\eta^*)\big)$. 
The map $\theta$ sends $\overline{f}$
to zero, hence we have  $\sum_{i=1}^\mu f_it_i = 0 \  {\rm modulo} \ 
\big(L_{\calO, \sigma} + I(\bbB_\calO)+(g^*_1, \dots, g_\nu^*) \big)$.
By what we have proved before, $\overline{\Phi}^{\mathstrut}$ 
is free with basis $\overline{\calO}$,
hence we deduce that $f_i \in   L_{\calO, \sigma} + I(\bbB_\calO)$, hence
$f_i \in \big(L_{\calO, \sigma} + I(\bbB_\calO) \big)\cap K[S_{c\calO,\sigma}]$
 for $i =1, \dots \mu$.
The equality 
$\big(L_{\calO, \sigma} + I(\bbB_\calO) \big)\cap 
K[S_{c\calO,\sigma}] = I(\bbG_{\calO, \sigma})$
yields the conclusion and the proof of b) is complete.

The proof of c) uses the same argument as above, 
which shows that if $\theta(\overline{f}) = 0$
then $\overline{f} = 0$.
Finally we prove d). At this point we know that
diagram $(1)$ is commutative, all the homomorphisms are 
injective, and both $\Psi$ and $\overline{\Phi}$ 
are free with basis $\overline{\calO}$.
Due to this particular structure, the surjectivity of ~$\theta$ 
is equivalent to the surjectivity of $\phi$.
Since  all the indeterminates which generate
$L_{\calO,\sigma}$ are killed,
we need to show that all the indeterminates
in $S_{\calO,\sigma}$ can be expressed as
polynomial functions of the indeterminates
in $S_{c\calO,\sigma}$.
We consider the generic
$(\calO, \sigma)$-Gr\"obner prebasis $\{g_1^*, \dots, g_\eta^*\}$ and
argue as in the proof of Proposition~\ref{iboHom}.
For every $j = \eta+1,\dots,\nu$ we produce elements
$b_j -\sum_{ \{ i \ | \ b_j >_{_\sigma} t_i \} }h_{ij}t_i $
which are in the ideal 
$(g_1^*, \dots, g_\eta^*) \subseteq (g_1^*, \dots, g_\nu^*)$.
Consequently, modulo 
$\big(L_{\calO, \sigma} + I(\bbB_{\calO})+(g^*_1,\dots, g^*_\nu )\big)$
we have $b_j -\sum_{ \{ i \ | \ b_j >_{_\sigma} t_i \} }h_{ij}t_i = 0$
as well as $b_j -\sum_{ \{ i \ | \ b_j >_{_\sigma} t_i \} }c_{ij}t_i = 0$
for every~$j = \eta+1,\dots,\nu$.
We deduce the relations $\sum_{ \{ i \ | \ b_j >_{_\sigma} t_i \} }
(c_{ij}-h_{ij})t_i = 0$ in
$U_{\calO} /\overline{L}_{\calO,\sigma}$
for every $j = \eta+1,\dots,\nu$. 
Using a)
we get the relations $c_{ij}=h_{ij}$
in the ring $B_{\calO}/\overline{L}_{\calO,\sigma}$,
for every $c_{ij} \in S_{\calO,\sigma}\setminus S_{c\calO,\sigma}$,
and every $j = \eta+1,\dots,\nu$, and the proof is complete.
\endproof

\medskip 
\begin{remark}\label{isoandiso}
After the  theorem, diagram $(1)$ can be rewritten in the following way.

$$
\begin{matrix}
G_{\calO,\sigma} &\iso{\phi}& B_{\calO}/\overline{L}_{\calO,\sigma}     \cr
&&& \cr 
\mapdown{\Psi} && \mapdown{\overline{\Phi}}  \cr
&&\cr
U_{\calO,\sigma}\ &\iso{\theta} &
U_{\calO}/\overline{L}_{\calO,\sigma} 
\end{matrix}
\eqno (2)
$$
\end{remark}

\goodbreak

\begin{corollary}\label{trueGfamily}
Let $\calO=\{t_1,\dots,t_\mu\}$  be an order ideal of monomials in $P$ and
let $\sigma$ be a term ordering on $\bbT^n$.

\begin{items}
\item The affine scheme $\bbG_{\calO, \sigma}$ parametrizes all 
zero-dimensional ideals $I$ in $P$ for which $\calO = \calO_\sigma(I)$.

\item The fibers over the $K$-rational points of the 
universal $(\calO, \sigma)$-Gr\"obner family
$\Psi: G_{\calO,\sigma} \To U_{\calO,\sigma}$
are the quotient rings $P/I$ for which $I$ is a zero-dimensional ideal
with the property that $\calO = \calO_\sigma(I)$. 
Moreover, the reduced $\sigma$-Gr\"obner basis of $I$ is obtained by 
specializing the $(\calO, \sigma)$-Gr\"obner prebasis $\{g^*_1, \dots, g^*_\eta\}$
to the corresponding maximal linear ideal.

\end{items}
\end{corollary}

\proof
The freeness of $\Psi$ implies that a) follows from b), and we prove b) in two steps.
A $K$-rational point of the universal $(\calO, \sigma)$-Gr\"obner basis 
family can be viewed as the $\Psi$-fiber over a maximal linear ideal
of $G_{\calO, \sigma}$. The latter is the 
 canonical projection of a maximal  linear ideal
$\goth{n}= (c_{ij}-a_{ij} \ | \ c_{ij} \in S_{c\calO, \sigma}, \ a_{ij} \in K)$
of $K[S_{c\calO, \sigma}]$. 
Let us put 
$IndL=\{(i,j) \ | \ c_{ij}\hbox{\ \rm  is a generator of\ } L_{\calO, \sigma} \}$. 
The theorem implies that 
then $\goth{n}$ is the contraction to~$K[S_{c\calO, \sigma}]$ 
of a maximal linear ideal 
$$
\goth{m}= \big(c_{ij}-a_{ij} \ | \  c_{ij} \in \{c_{11}, \dots, c_{\mu\nu}\}, 
\ a_{ij} \in K, \ a_{ij}= 0 {\rm \ for\  all \ } (i,j) \in IndL \big)
$$ 
of $K[c_{11}, \dots, c_{\mu\nu}]$.
The ideal $\goth{m}$ contains $I(\bbB_\calO)$, hence if we substitute
$c_{ij}$ with~$a_{ij}$  in the polynomials $g_1^*, \dots, g_\nu^*$, we get 
polynomials $\overline{g}_1, \dots, \overline{g}_\nu$ in $P$ which form
the~$\calO$-border basis of the ideal 
$I  = (\overline{g}_1, \dots, \overline{g}_\nu)$. 
Moreover, by construction
we have $\LT_\sigma(\overline{g_j}) = b_j$ for $j = 1,\dots, \nu$.
Hence $\{\overline{g}_1, \dots, \overline{g}_\eta\}$ is the 
reduced $\sigma$-Gr\"obner basis of $I$ by Proposition 6.4.18 of ~\cite{KR2}.

Conversely, let $I$ be a zero-dimensional ideal in $P$ such 
that $\calO_\sigma(I) = \calO$
and let $\{g_1, \dots, g_\eta\}$ be its reduced $\sigma$-Gr\"obner basis. 
Using the division algorithm,
we represent all the elements in $\partial\calO \setminus c\calO$ 
uniquely (modulo $I$)  as linear combinations 
of elements in $\calO$. In this way, the $\calO$-border 
basis  $(g_1, \dots, g_\nu)$ of $I$ is constructed.
Collecting the coefficients, we produce a maximal linear ideal in $B_\calO$, 
equivalently a  rational point~$\mathbf p$ of  $\bbB_\calO$. 
By construction, $b_j = \LT_\sigma(g_j)$ for $j = 1, ..., \nu$, 
and hence the coordinates of $\mathbf p$ corresponding to
the indices $ij$ such that  $(i,j) \in IndL$ have to be zero.
In conclusion, the point $\mathbf p$ corresponds to a maximal linear
ideal $\goth{m}$ of $B_\calO/{\overline{L}_{\calO, \sigma}}$ hence to a 
maximal linear ideal of $G_{\calO, \sigma}$ by the theorem, hence to a rational 
point~$\mathbf q$ of~$\bbG_{\calO, \sigma}$.
The ideal itself is represented via its reduced $\sigma$-Gr\"obner basis 
$\{\overline{g}_1, \dots, \overline{g}_\eta\}$ in  the 
$\Psi$-fiber over~$\goth{m}$.
\endproof

\begin{remark}\label{specandproj}
Diagram $(2)$ gives rise to the corresponding 
diagram 
$$
\begin{matrix}
\bbG_{\calO,\sigma} &\cong& \Spec(B_\calO/{\overline{L}_{\calO, \sigma}})    \cr
&&& \cr 
\mapup{\pi_\Psi} && \mapup{\pi_{\overline{\Phi}}}  \cr
&&\cr
\Spec(U_{\calO,\sigma})\ &\cong &
\Spec(U_\calO/{\overline{L}_{\calO, \sigma}}) 
\end{matrix}
\eqno (3)
$$
of affine schemes, but more can be said.
Let $W$, $\overline{W}$, $V$ be systems of positive weights,
chosen suitably to  satisfy Theorem~\ref{homogeneousfamily}.
Then $G_{\calO, \sigma}$ is a $W$-graded ring, $B_\calO$ is a 
 $\overline{W}$-graded ring, $U_{\calO, \sigma}$ is a $(V,W)$-graded ring,
 and $U_\calO/{\overline{L}_{\calO, \sigma}}$ is a $(V, \overline{W})$-graded ring.

With the above assumptions  we see that 
diagram $(2)$ gives rise to a 
diagram
$$
\begin{matrix}
\Proj(G_{\calO, \sigma})&\cong& \Proj(B_\calO/{\overline{L}_{\calO, \sigma}})    \cr
&&& \cr 
\mapup{\Pi_\Psi} && \mapup{\Pi_{\overline{\Phi}}}  \cr
&&\cr
\Proj(U_{\calO,\sigma})\ &\cong &
\Proj(U_\calO/{\overline{L}_{\calO, \sigma}}) 
\end{matrix}
\eqno (3)
$$

\goodbreak

\noindent of projective schemes 
$\Proj(G_{\calO, \sigma})$,  
$\Proj(B_\calO/{\overline{L}_{\calO, \sigma}})$,  
$\Proj(U_{\calO,\sigma})$, 
$\Proj(U_\calO/{\overline{L}_{\calO, \sigma}})$
such that 
${\Proj(G_{\calO, \sigma})\subset \bbP(W)}$,  
$\Proj(B_\calO/{\overline{L}_{\calO, \sigma}})   \subset \bbP(\overline{W})$,  
$\Proj(U_{\calO,\sigma}) \subset \bbP(V,W)$, and 
$\Proj(U_\calO/{\overline{L}_{\calO, \sigma}})  \subset \bbP(V, \overline{W})$
where $\bbP(W)$, $\bbP(\overline{W})$, $\bbP(V,W)$, and $\bbP(V, \overline{W})$ 
are the corresponding  weighted projective spaces.
 
Moreover, let $\mathbf p = (a_{ij}) \in \bbG_{\calO, \sigma}$ be a rational point, 
let $I\subset P$ be the corresponding ideal 
according to Corollary~\ref{trueGfamily}, 
let $v_i = \deg(x_i)$ in the $V$-grading, and let $w_{ij} = \deg(c_{ij})$ in the $W$-grading. 
Then it is well-known that the substitution 
$a_{ij} \To t^{w_{ij}}a_{ij}$ gives rise to a  flat family of ideals
whose 
general fibers are ideals isomorphic to $I$, and whose special 
 fiber is the monomial ideal $\LT_\sigma(I)$. In the setting of diagram $(2)$, the 
 rational monomial curve which parametrizes such family is a 
 curve in $\bbG_{\calO, \sigma}$ which connects the two points 
 representing $I$ and $\LT_\sigma(I)$. In the setting of diagram $(3)$, 
 the rational monomial curve is simply a point 
 in ${\Proj(G_{\calO, \sigma})\subset \bbP(W)}$, which represents all the above ideals 
 except the special one.

\end{remark}

\goodbreak

\section{Consequences and problems}
\label{Applications}

We open the section by discussing the relation between our 
construction of~$I(\bbG_\calO)$ and other constructions 
described in the literature
(see for instance~\cite{CV} and~\cite{NS}). If one starts with the
generic $\sigma$-Gr\"obner prebasis $\{g^*_1, \dots, g^*_\eta\}$ one can 
construct an affine subscheme of $\bbA^{s(c\calO, \sigma)}$ 
in the following way.
Using Buchberger Algorithm one reduces  the critical pairs 
of the leading terms of the $\sigma$-Gr\"obner prebasis as much as possible.
The reduction stops when a polynomial is obtained which is a linear combination
of the elements in $\calO$ with coefficients in $K[S_{c\calO, \sigma}]$.
Collecting all coefficients obtained in this way for all the critical pairs, 
one gets a set which generates an ideal $J$ in $K[S_{c\calO, \sigma}]$. 
Clearly each zero of $J$ gives rise to a specialization of the generic 
$\sigma$-Gr\"obner prebasis which is, by construction, the 
reduced $\sigma$-Gr\"obner basis of a zero-dimensional 
ideal $I$ in $P$ for which $\calO = \calO_\sigma(I)$.
However, there is a drawback; the reduction procedure in Buchberger Algorithm 
is far from being unique. This observation leads to the following definition which
puts the above description in a more formal context.

\begin{definition}\label{reductionideal}
Let $J$ be an ideal in  $K[S_{c\calO, \sigma}]$  such that 
$\Spec(K[S_{c\calO, \sigma}]/J)$ 
param\-etrizes all zero-dimensional ideals~$I$ in $P$ 
for which $\calO = \calO_\sigma(I)$. Then $J$ is called 
an {\it $(\calO, \sigma)$-parametrizing ideal}.

Let $J$ be an an  $(\calO, \sigma)$-parametrizing ideal and
assume that  there exists a finite set $\calS$ of polynomials of  type 
$\sum_{j=1}^\eta f_jg^*_j = \sum_{t_i \in \calO} r_i t_i  $
where the $f_i$'s are polynomials in $K[\mathbf x, S_{c\calO, \sigma}]$, the 
$r_i$'s are polynomials in $K[S_{c\calO, \sigma}]$, and $J$ is 
generated by the~$r_i$'s. 
Then $J$ will be called an {\it $(\calO, \sigma)$-reduction ideal}, 
and $\calS$ an {\it $(\calO, \sigma)$-reduction set}\/ of~$J$.
\end{definition}

\begin{lemma}\label{reduction}
Let $J$ be an $(\calO, \sigma)$-parametrizing ideal. Then 
the canonical homomorphism
$
K[S_{c\calO, \sigma}]/J \To   
K[\mathbf x, S_{c\calO, \sigma}]/\big( J + (g^*_1, \dots, g^*_\eta)\big)
$
makes the quotient ring 
$K[\mathbf x, S_{c\calO, \sigma}]/\big( J + (g^*_1, \dots, g^*_\eta)\big)$ 
into a free $K[S_{c\calO, \sigma}]/J$-module, and  
a basis is the set of  the residue classes of the elements of $\calO$.

\end{lemma}

\proof
To prove this lemma we have to show that the residue clases of 
the elements in $\calO$ generate
$K[\mathbf x, S_{c\calO, \sigma}]/\big( J + (g^*_1, \dots, g^*_\eta)\big)$ 
and are linearly independent.

{\it $\overline{\calO}$ generates}. It is enough to use the 
$\overline{\sigma}$-division algorithm, as we did in the proof
of Theorem~\ref{homogeneousfamily}.

{\it $\overline{\calO}$ is  linearly independent over 
$K[S_{c\calO, \sigma}]/J$}.  Suppose not.  
Then there would be a non-empty  open set of 
$\Spec(K[S_{c\calO, \sigma}]/J)$, whose 
maximal linear ideals would represent 
ideals~$I$ of $P$ for which 
$\calO_\sigma(I) \subset \calO$,  a contradiction.
\endproof

\begin{remark}
If $J$ is an  $(\calO, \sigma)$-parametrizing ideal, then
it is not necessarily an~$(\calO, \sigma)$-reduction ideal. 
It suffices to pick an ideal $J$ which is an  
$(\calO, \sigma)$-param\-etrizing ideal in  
$K[S_{c\calO, \sigma}]$ but not radical. 
Then $\sqrt{J}$ is still an $(\calO, \sigma)$-parametrizing 
ideal but not necessarily an~$(\calO, \sigma)$-reduction ideal. 
\end{remark}

\begin{lemma}\label{substelim}
Let $\mathbf{x} = x_1, \dots, x_n$, $\mathbf{y} = y_1, \dots, y_m$, and let
$P = K[\mathbf x]$, $Q = K[\mathbf x, \mathbf y]$.
Let $g_1, \dots, g_t $ be polynomials in $Q$, let $J$ be 
the ideal generated by $\{g_1, \dots, g_t\}$, and assume 
that there exist polynomials $f_1(\mathbf x), \dots, f_m(\mathbf x)$ such that 
the elements ${y_1-f_1(\mathbf x), \dots, y_m - f_m(\mathbf x)}$ are in $J$. 
Then the ideal  $J\cap K[\mathbf x]$ is generated by 
$\{g_1(\mathbf x, \mathbf f), \dots, g_t(\mathbf x, \mathbf f) \}$
where $\mathbf f = (f_1, \dots, f_m)$.
\end{lemma}

\proof
Every polynomial $g \in Q$ can be written as
$$
g(\mathbf x, \mathbf y) = \sum_{i=1}^mh_i(y_i-f_i) + g(\mathbf x, \mathbf f)
$$
and the remainder $g(\mathbf x, \mathbf f)$ is unique, since
$\{y_1-f_1(\mathbf x), \dots, y_m - f_m(\mathbf x)\}$ is a Gr\"obner basis with respect to an ordering which eliminates $\mathbf y$. Now the conclusion follows easily.
\endproof

\goodbreak

\begin{proposition}\label{BBisreduction}
The ideal $I(\bbG_{\calO, \sigma})$ is an $(\calO, \sigma)$-reduction ideal. 
\end{proposition}

\proof Following Definition~\ref{reductionideal} we have to prove that 
$I(\bbG_{\calO, \sigma})$ is 
an $(\calO, \sigma)$-param\-etrizing ideal,  and that
there exists an $(\calO, \sigma)$-reduction set of~$I(\bbG_{\calO, \sigma})$.
The first claim was proved in Corollary~\ref{trueGfamily}. To prove the second claim we 
use~\cite{KR3}, Proposition 4.1 and~\cite{KK1},  Section 4 to get 
generators 
of the ideal $I(\bbB_\calO)$ as the coefficients of the $t_i$'s in polynomial expressions
of type $\sum_{j=1}^\nu f_jg_j = \sum_{t_i \in \calO} r_i t_i $ where
the $f_i$'s are polynomials in $K[\mathbf x, \mathbf c]$ and  the 
$r_i$'s are polynomials in $K[\mathbf c]$.  

Consequently, to get generators
of  the ideal $L_{\calO, \sigma}+I(\bbB_\calO)$
we pick all these polynomial expressions,  
set equal to zero all the indeterminates
which generate $L_{\calO, \sigma}\!$, and get expressions  
$\sum_{j=1}^\nu f^*_jg^*_j = \sum_{t_i \in \calO} r^*_i t_i $
where the $f^*_i$'s are polynomials in $K[\mathbf x, S_{\calO, \sigma}]$ and  the 
$r^*_i$'s are polynomials in $K[S_{\calO, \sigma}]$.  

For the sake of clarity,  let us call~$\mathbf{\tilde{c}}$ 
the set of indeterminates $S_{c\calO, \sigma}$, and 
$\mathbf{\tilde{d}}$ the set of indeterminates 
$S_{\calO, \sigma}\setminus S_{c\calO, \sigma}$.
We rewrite $\sum_{j=1}^\nu f^*_jg^*_j = \sum_{t_i \in \calO} r^*_i t_i $ 
$$
\sum_{j=1}^\eta f^*_j(\mathbf{\tilde{c}}, \mathbf{\tilde{d}})\, g^*_j(\mathbf{\tilde{c}})+ 
\sum_{j=\eta+1}^\nu f^*_j(\mathbf{\tilde{c}}, \mathbf{\tilde{d}})\, g^*_j(\mathbf{\tilde{d}}) =
\sum_{t_i \in \calO} r^*_i(\mathbf{\tilde{c}}, \mathbf{\tilde{d}})\, t_i   \eqno(*)
$$
Once more, we argue as in the proof of Theorem~\ref{GandB}.d, and 
for every $j = \eta+1,\dots,\nu$ we produce elements
$b_j -\sum_{ \{ i \ | \ b_j >_{_\sigma} t_i \} }h_{ij}t_i$
which are in the ideal 
$(g_1^*, \dots, g_\eta^*)$.
Using Theorem~\ref{GandB}, we get  relations $c_{ij} = h_{ij}$ in the ring
$K[\mathbf c]/(L_{\calO, \sigma}+I(\bbB_\calO)$, 
in other words  relations $c_{ij}-h_{ij} \in L_{\calO, \sigma}+I(\bbB_\calO)$.
We make the substitution $c_{ij} \To h_{ij}$ in the 
expressions $(*)$, write $\mathbf{h_{ij}}$ for the tuple
$h_{ij}$, and get expressions
$$
\sum_{j=1}^\eta \tilde{f}_j(\mathbf{\tilde{c}}) g^*_j(\mathbf{\tilde{c}}) =
\sum_{t_i \in \calO} \tilde{r}_i(\mathbf{\tilde{c}}) t_i  
\hbox{\quad where \quad} \tilde{r}_i(\mathbf{\tilde{c}}) = r^*_i(\mathbf{\tilde{c}}, \mathbf{h_{ij}})
 \eqno(**)
$$
Now it suffices to prove that the set of all the $ \tilde{r}_i(\mathbf{\tilde{c}})$ generates $I(\bbG_{\calO,\sigma})$.
We recall the equality $I(\bbG_{\calO,\sigma}) = 
\big(L_{\calO, \sigma}+I(\bbB_\calO) \big)\cap K[S_{c\calO, \sigma}]$ 
and we know that the ideal $L_{\calO, \sigma}+I(\bbB_\calO)$ is generated by $L_{\calO, \sigma}$
and the polynomials $r^*_i(\mathbf{\tilde{c}}, \mathbf{\tilde{d}})$, hence the conclusion follows from the lemma.
\endproof
%


\begin{proposition}\label{democracy}
All the $(\calO, \sigma)$-reduction ideals are equal. 
\end{proposition}

\proof
Let $J_1$, $J_2$ be $(\calO, \sigma)$-reduction ideals. 
By interchanging the role of $J_1$ and~$J_2$ it suffices to 
prove that $J_1\subseteq J_2$.  Let $\calS$ 
be an $(\calO, \sigma)$-reduction set  of~$J_1$. Every element in $\calS$ has the shape
$\sum_{j=1}^\eta f_jg^*_j = \sum_{t_i \in \calO} h_i t_i$.
We consider the canonical homomorphism 
$$
K[S_{c\calO, \sigma}]/J_2 \To   
K[\mathbf x, S_{c\calO, \sigma}]/\big( J_2 + (g^*_1, \dots, g^*_\eta)\big)
$$
and deduce that $\sum_{t_i \in \calO} h_i t_i =0$ in the ring
$K[\mathbf x, S_{c\calO, \sigma}]/\big( J_2 + (g^*_1, \dots, g^*_\eta)\big)$ which is free
over $K[S_{c\calO, \sigma}]/J_2 $ by Lemma~\ref{reduction}. Therefore the coefficients $h_i$
are zero in  the ring $K[S_{c\calO, \sigma}]/J_2$. In particular, they belong to
$J_2$ and the proof is complete.
\endproof

\goodbreak

A combination of Theorems~\ref{homogeneousfamily} and~\ref{GandB}
yields a remarkable property of   ${\mathbb G}_{\calO, \sigma}$.
A~similar result can be found in~\cite{NS} proposition 4.3. 
The main difference is that there the authors deal 
with standard homogeneous saturated ideals. 
Moreover their proof is incorrect.

\begin{corollary}\label{affinecell}
Let $\calO \subset \bbT^n$ be an order ideal of monomials,
let $\sigma$ be a term ordering on $\mathbb{T}^n$,
and let $\mathbf o$ be the origin in the affine space 
$\bbA^{s(c\calO,\sigma)}$.

\begin{items}
\item
The point $\mathbf o$ belongs to ${\mathbb G}_{\calO, \sigma}$.

\item The following conditions are equivalent
\begin{items}
\item[1)] The scheme ${\mathbb G}_{\calO, \sigma}$ is 
isomorphic to an affine space.
\item[2)]  The point $\mathbf o$
is a smooth point of ${\mathbb G}_{\calO, \sigma}$.

\end{items}

\end{items}
\end{corollary}

\proof
The point $\mathbf o$ corresponds to the monomial ideal generated by $c\calO$,
and hence it belongs to ${\mathbb G}_{\calO, \sigma}$ by Corollary~\ref{trueGfamily}.
To prove part b), it is clearly sufficient to show that $2)$ implies $1)$. We argue as follows.
Suppose that among the $W$-homogeneous generators of the ideal $I(\bbG_{\calO, \sigma})$ 
there is one, say $f$, of  type $c_{ij} - g$ with the property  that $c_{ij}$ does not divide 
any elements in the support of $g$.  The graded ring $G_{\calO, \sigma}/(f)$ is isomorphic
to a graded $K$-algebra embedded  in a polynomial ring with one less indeterminate, 
the isomorphism being constructed by substituting~$c_{ij} $ with $g$.
Suppose we do this operation until no polynomial like $f$ is found anymore,  
call $Q/J$ the graded algebra obtained in this way, with $Q$ a polynomial ring, 
and~$J$ a homogeneous ideal.
We claim that no polynomial in $J$ can have a non-zero linear part. 
For contradiction, suppose that a polynomial $h$ of that type exists,
and let $c_{ij}$  be an indeterminate in the support of the linear part of~$h$. Then~$c_{ij}$ must divide
another power product in the support of $h$ which is impossible since~$J$ is homogeneous
with respect to a set of positive weights. In conclusion, we have~$J = (0)$.
 \endproof
 
The  algebraic argument given in the above proof agrees with the well-known fact that
a quasi-cone over a projective subscheme $\mathbb X$ of a weighted projective 
scheme $\bbP(V)$ is smooth if and only if $\mathbb X = \bbP(V)$.

\begin{remark}
There is a  strong difference between 
$\bbG_{\calO, \sigma}$ 
and $\mathbb {B}_{\calO, \sigma}$ even when $n=2$.
It is known that for $n=2$ the scheme 
$\mathbb {B}_{\calO}$ is smooth and irreducible.
However, unlike the case of~$\bbG_{\calO, \sigma}$ as
explained in Corollary~\ref{affinecell}, it does not need to be an affine cell
(i.e.\ isomorphic to an affine space)
as the following example shows.
\end{remark}

\goodbreak

\begin{example}\label{square}
This is an example  where~$\bbG_{\calO, \sigma}$ 
is isomorphic to an affine space of 
dimension~$9$, and where $\bbB_{\calO, \sigma}$ 
is a smooth irreducible variety
of dimension $10$ not isomorphic to an affine space. 
Let $P=k[x, y]$ and 
$\calO = (1,\, x, \, y, \, x^2, \, y^2)$. 
Then~$\partial{\calO} = (xy,\, y^3, \, x^3, \, xy^2,\, x^2y)$ and so 
$\mu = \nu = 5$.  
Using \cocoa\ we compute~$I(\bbB_\calO)$ and find out that 
$\dim(B_\calO) = 10$. It is the expected number
since the Hilbert scheme has only one component whose 
general point corresponds to the ideal of five distinct points 
in $\bbA^2$, and hence depends on ten parameters.
Moreover we see that~${\mathbb B}_\calO$ is isomorphic to
a smooth irreducible variety
of dimension $10$, embedded in an affine space of dimension 
$14$ and described by an ideal with $9$ generators. 
Looking at the shape of the equations it is 
easy to see that it is not isomorphic to an affine space.
The fact that $B_\calO$ is smooth and irreducible agrees with 
a general statement that all the border basis schemes in 
two indeterminates are smooth and irreducible 
(see~\cite{Ha} Proposition 2.4 and~\cite{Hu3} Corollary 9.5.1).

Now we let $\sigma = {\tt DegLex}$. 
We see that 
$S_{\calO, \sigma} = \{c_{11}, \dots c_{55}\}\setminus \{c_{41}\}$
since ${xy <_\sigma x^2}$,
hence $L_{\calO, \sigma}= (c_{41})$.
Then  we check that 
$c\calO = (xy,\,  y^3, \, x^3)$, hence $\eta = 3$, and~$L_{}$ 
is the ideal generated by $c_{41}$ in $K[c_{11}, \dots, c_{53}]$.
Now we check with \cocoa\  that the ring $B_\calO/\overline{L}_{\calO, \sigma}$
is isomorphic to a polynomial ring with $9$ indeterminates, 
and, in agreement with Corollary~\ref{trueGfamily}, we deduce from the 
Theorem~\ref{GandB} that also
$G_{\calO, \sigma}$ is isomorphic to a 
polynomial ring with $9$ indeterminates.
\end{example}

As a natural follow up to Theorem~\ref{GandB} we look for conditions under which
we have 
$\overline{L}_{\calO, \sigma} = \emptyset$ so that  diagram $(2)$
can be written as
$$
\begin{matrix}
G_{\calO,\sigma} &\iso{\phi}& B_{\calO}   \cr
&&& \cr 
\mapdown{\Psi} && \mapdown{\overline{\Phi}}  \cr
&&\cr
U_{\calO,\sigma}\ &\iso{\theta} &
U_{\calO}
\end{matrix}
\eqno (4)
$$
In other words we  look for conditions under 
which the border basis scheme and the Gr\"obner basis scheme are isomorphic.
We recall some definitions from~\cite{KR3} (Definition 2.7) and~\cite{OS}.

\begin{definition}
Let $\calO$ be an order ideal, let $V$ be a matrix in $\Mat_{1,n}(\mathbb N_+)$,
and let $\sigma$ be a term ordering on $\bbT^n$.

\begin{items}
\item The order ideal
$\calO$ is said to have a $maxdeg_V\ border$ if
$\deg_V(b) \ge \deg_V(t)$
for every $b \in c\calO$ and every $t \in \calO$.

\item Similarly, $\calO$ is said to be 
a {\it  $V$-cornercut}\/  (or to have a $strong\ maxdeg_V\ border$) if
$\deg_V(b) > \deg_V(t)$
for every $b \in c\calO$ and every $t \in \calO$.

\item 
The order ideal  $\calO$ is said to be a {\it  $\sigma$-cornercut}\/ if
$b>_\sigma t$ for every $b \in c\calO$ and every $t \in \calO$.

\end{items}
\end{definition}

\begin{proposition}\label{iso}
Let $\calO$ be an order ideal and $\sigma$ a term 
ordering on $\bbT^n\!$. Consider the following conditions.
\begin{items}
\item[$a_1)$] The 
canonical embedding of $K[S_{\calO,\sigma}]$ in 
$K[c_{11}, \dots, c_{\mu\nu}]$ induces an 
isomorphism between $G_{\calO, \sigma}$ and $B_\calO$.
\item[$a_2)$]  The 
canonical embedding of $K[\mathbf x, S_{\calO,\sigma}]$ in 
$K[\mathbf x, c_{11}, \dots, c_{\mu\nu}]$ induces an 
isomorphism between $U_{\calO, \sigma}$ and $U_\calO$.
\item[$b_1)$]  The ideal $L_{\calO, \sigma}$ is the zero ideal.
\item[$b_2)$]  The order ideal $\calO$ is a $\sigma$-cornercut.
\end{items}
Then $a_1)$ is equivalent to $a_2)$, $b_1)$ is equivalent to $b_2)$, and
$b_1)$ implies $a_1)$.
\end{proposition}
\proof
The equivalence of $a_1)$ and $a_2)$ follows from 
Theorem~\ref{GandB}, since $U_\calO$ is a free~$B_\calO$ module with 
basis $\overline{\calO}$, and also $U_{\calO,\sigma}$ is 
a free $G_{\calO,\sigma}$ module with 
basis~$\overline{\calO}$. Next we prove the implication $b_1)\implies b_2)$. 
If $L_{\calO, \sigma}$ is the zero ideal, then~$b_j >_\sigma t_i$ 
for every $j = 1, \dots, \nu$ and every $i = 1, \dots, \mu$. Consequently
we have $b>_\sigma t$ for every $b \in c\calO$ and 
every $t \in \calO$ i.e.\ $\calO$ is a $\sigma$-cornercut. 
The implication $b_2) \implies b_1)$ follows from the definition 
of $L_{\calO, \sigma}$ and the
implication $b_1)\implies a_1)$ follows immediately from Theorem~\ref{GandB}.
\endproof

\begin{remark}
Let us make some remarks about this proposition.
\begin{items}
\item Remark~\ref{posweights} has the following implication.
If condition $b_2)$ is 
fulfilled  i.e.\ $\calO$ is a $\sigma$-cornercut, 
then there exists a system $V$ of positive weights such that
$\calO$ is a $V$-cornercut.

\item Example~\ref{square} shows that in the above proposition one cannot 
substitute condition $b_2)$ with the weaker condition that  the order
ideal $\calO$ has a  maxdeg$_V$-border. 
 
\item The author does not know whether all the conditions 
of the above proposition are equivalent.
\end{items}
\end{remark}

As a consequence of Proposition~\ref{iso} 
we give a very short proof of the fact that if~$\calO$ has the shape of a segment
then $B_\calO$ is an affine space.

\begin{corollary}\label{stick}
Let $\calO = \{1, x_n, x_n^2, \dots, x_n^{\mu-1}\} \subset \bbT^n$. 
Then $\bbB_\calO$ is isomorphic to the affine space $\bbA^{\mu n}$.
\end{corollary}

\proof
Clearly $\calO$ is a ${\tt Lex}$-cornercut, hence $B_\calO$ is isomorphic to 
$G_{\calO, {\tt Lex}}$ and we have $g^*_j = g_j$ for $j = 1,\dots, \nu$. 
Corollary~\ref{trueGfamily} implies that $G_{\calO, {\tt Lex}}$
parametrizes all zero-dimensional ideals~$I$ in $P$ for which $\calO =\calO_{\tt Lex}(I)$.
Hence $I(\bbG_{\calO,\tt Lex})$ contains relations under 
which the generic $\tt Lex$-Gr\"obner prebasis
is the reduced~$\tt Lex$-Gr\"obner basis of an
ideal $I$ in $P$ for which $\calO =\calO_{\tt Lex}(I)$. 
On the other hand,  it is clear that $\eta = n$ and the  
generic $\tt Lex$-Gr\"obner prebasis 
consists of  $n$ polynomials whose leading terms are
$x_1, \dots, x_{n-1}, x_n^{\mu}$. 
They are pairwise coprime, hence  every specialization
of the generic $\tt Lex$-Gr\"obner prebasis is a reduced
$\tt Lex$-Gr\"obner basis.
It follows that $I(\bbG_\calO)$ is the zero
ideal and the proof is complete.
\endproof

We observe that the explicit isomorphism of $B_\calO$ with the polynomial ring 
$K [x_{11}, \dots, x_{\mu n}]$ is given by expressing  the indeterminates  
$x_{1, n+1},\dots, x_{\mu\nu}$ as
polynomials in the indeterminates in $x_{11}, \dots, x_{\mu n}$,
as explained in the proof of Theorem~\ref{GandB}.d.

\bigskip

The final part of the section and hence of the paper 
is devoted to a general remark and the discussion of some open problems.

\begin{remark} In the paper~\cite{KR3} we have introduced 
and discussed the {\it homogeneous border basis scheme}. 
With the obvious modifications one can as well introduce 
the {\it homogeneous Gr\"obner basis scheme}. 
\end{remark}

 Using Theorem~\ref{GandB}
and Remark~\ref{isoandiso} we know the precise 
relation between the two schemes  $\bbG_{\calO, \sigma}$ and $\bbB_\calO$.
It is then quite natural to ask the following question.

\medskip

\noindent {\bf Question 1:}{\it \ 
Is there any connection between the smoothness of the origin 
in~$\bbG_{\calO, \sigma}$ and the smoothness of the origin in $\bbB_\calO$?}

\medskip

The scheme $\bbG_{\calO, \sigma}$ is connected since
it is a quasi-cone,  and hence all its points are connected to 
the origin  (see Remark~\ref{specandproj}). 
However, the problem of the connectedness of $\bbB_\calO$ is still open, so let me state it formally. 

\medskip

\noindent {\bf Question 2:}{\it \ 
Is $\bbB_\calO$ connected?}

\bigskip

\subsection*{Acknowledgements}
This work was almost entirely conducted while I was visiting
the Indian Institute of Science in Bangalore, India. 
It is my pleasure to thank the institution, in particular Professor Dilip Patil for the 
warm hospitality and  for providing an exceptionally good atmosphere.

\bigbreak


\begin{thebibliography}{99}


\bibitem{Bay} D.\ Bayer, The division algorithm and the Hilbert scheme.
Thesis, Harvard  University, Cambridge, MA, (1982).


\bibitem{CoCoA} The \cocoa\ Team, {\it \cocoa : a system for doing
Computations in Commutative Algebra},
available at {\tt http://cocoa.dima.unige.it}.

\bibitem{CV} A.\ Conca and G.\ Valla,
Canonical Hilbert-Burch matrices for ideals of $k[x,y]$, \\
{\tt arXiv:math$\!\backslash\!$0708.3576}.

\bibitem{E} D.\ Eisenbud, {\it Commutative Algebra with a View Toward
Algebraic Geometry}, Springer, Heidelberg 1995.

\bibitem{GLS} T.S.\ Gustavsen, D.\ Laksov and R.M.\ Skjelnes,
An elementary, explicit proof of the existence of Hilbert
schemes of points,  {\tt arXiv:math$\!\backslash\!$0506.161v1}.

\bibitem{Ha} M.\ Haiman, q,t-Catalan numbers and the Hilbert
scheme, Discr.\ Math.\ {\bf 193} (1998), 201--224.


\bibitem{Hu1} M.\ Huibregtse, A description of certain affine open
schemes that form an open covering of~${\rm Hilb}^n_{\mathbb A^2_k}$,
Pacific J. Math. {\bf 204} (2002), 97--143.

\bibitem{Hu2} M.\ Huibregtse, An elementary construction of the
multigraded Hilbert scheme of points,
Pacific J. Math. {\bf 223} (2006), 269--315.

\bibitem{Hu3} M.\ Huibregtse, The cotangent space at a monomial 
ideal of the Hilbert scheme of points of an affine space
{\tt arXiv:math$\!\backslash\!$0506575}.



\bibitem{KK1} A. Kehrein and M. Kreuzer, Characterizations of border
bases,
J.\ Pure Appl.\ Alg.\ {\bf 196} (2005), 251--270.



\bibitem{KR1} M.\ Kreuzer and L.\ Robbiano, {\it Computational
Commutative Algebra 1},
Springer, Heidelberg 2000.

\bibitem{KR2} M.\ Kreuzer and L.\ Robbiano, {\it Computational
Commutative Algebra 2},
Springer, Heidelberg 2005.

\bibitem{KR3} M.\ Kreuzer and L.\ Robbiano, Deformations of border
bases, {\tt arXiv:0710.2641}, To appear in Collectanea Mathematica.


\bibitem{MS} E.\ Miller and B.\ Sturmfels, {\it Combinatorial
Commutative Algebra}, Springer, New York 2005.

\bibitem{M} B.\ Mourrain, A new criterion for normal form algorithms,
AAECC Lecture Notes in Computer Science \ {\bf 1719} (1999), 430--443.

\bibitem{NS} R.\ Notari and M. L.\ Spreafico, A stratification of
Hilbert schemes by initial ideals and applications,
Manuscripta Math.\ {\bf 101} (2000), 429--448.

\bibitem{OS} S.\ Onn and B.\ Sturmfels, Cutting corners,
Adv. Appl. Math.\ {\bf 23(1)} (1999), 29--48.

\bibitem{S} H.J.\ Stetter, {\it Numerical Polynomial Algebra}, SIAM,
Philadephia 2004.

\end{thebibliography}
\end{document}